\documentclass[10pt]{article}
\usepackage{amsmath,amsfonts}

\newtheorem{Theorem}{Theorem}[section]
\newtheorem{Definition}[Theorem]{Definition}
\newtheorem{Proposition}[Theorem]{Proposition}
\newtheorem{Lemma}[Theorem]{Lemma}

\newtheorem{Remark}[Theorem]{Remark}
\newtheorem{Example}[Theorem]{Example}

\makeatletter
\newif\ifmsbmloaded@

\@addtoreset{equation}{section}

\makeatother

\def \Espace{\renewcommand{\arraystretch}{1.7} }
\def\R{\mathbb R}
\def\N{\mathbb N}

\def\E{\mathbb E}
\def\P{\mathbb P}

\def\Lr{\mathcal L}

\def\ds{\displaystyle}
\def\vep{\varepsilon}

\begin{document}

 \begin{center} {\Large \bf Markov solutions for  the $3D$ stochastic
Navier--Stokes equations with state dependent noise}

\vspace{2mm}

Arnaud Debussche, Cyril Odasso

 IRMAR and  Ecole Normale Sup\'erieure de Cachan, antenne de
Bretagne, Campus de Ker Lann, 35170 Bruz, France
\end{center}\vspace{5mm}

\noindent{\bf Abstract: } We construct a Markov family of
solutions for the 3D Navier-Stokes equation perturbed by a non
degenerate noise. We improve the result of \cite{DPD-NS3D} in two
directions. We see that in fact not only a transition semigroup
but a Markov family of solutions can be constructed. Moreover, we
consider a state dependant noise. Another feature of this work is
that we  greatly simplify the proofs of  \cite{DPD-NS3D}.

\noindent {\bf 2000 Mathematics Subject Classification AMS}: 35Q10,
60H15,  37L40.


\section{Introduction}

We are concerned with the stochastic Navier--Stokes equations
on  $\mathcal O$ an open bounded domain of $\R^3$ with smooth
boundary $\partial \mathcal O.$ The unknowns are the velocity $X(t,\xi)$ and the pressure
 $p(t,\xi)$ defined for $t>0$ and $\xi\in \overline{\mathcal O}$:
\begin{equation}
\label{e1.1} \left\{\begin{array}{l} dX(t,\xi)=[\Delta
X(t,\xi)-(X(t,\xi)\cdot\nabla)X(t,\xi)]dt-\nabla
p(t,\xi)dt+f(\xi)dt+\Phi(X(t,\cdot))(\xi)\;dW,
\\
\mbox{\rm div}\;X(t,\xi)=0,
\end{array}\right.
\end{equation}
with Dirichlet boundary conditions
$$
X(t,\xi)=0,\quad t>0,\;\xi\in \partial \mathcal O,
$$
and supplemented with the initial condition
$$
X(0,\xi)=x(\xi),\;\xi\in  \mathcal O.
$$
We have taken the viscosity equal to $1$ since it plays no particular role in this work. The noise
is of white noise type and its covariance may depend on the noise through the operator
$\Phi(X)$.

It is classical that this equation, both in the deterministic and stochastic case, has a global
weak solution. Uniqueness of global solutions is a long standing open problem and one of the
main challenge of the theory of partial differential equations (see \cite{temam-survey} for
a survey on these questions). Although an intense research has been performed on this aspect
and many new ideas have appeared, this problem seems to be still out of reach.

In the stochastic case, a less difficult problem is the uniqueness
in law in the spirit of Stroock and Varadhan \cite{SV}. Some
progress have recently been obtained on that aspect. In
\cite{DPD-NS3D}, it has been proved that if the noise is additive
and sufficiently non degenerate it is possible to construct
directly a solution of the associate Kolmogorov equation.
Unfortunately, this solution is not sufficiently smooth and Ito
formula cannot be applied. Thus the Stroock and Varadhan program
cannot be accomplished and uniqueness in law remains open.
However, using other ideas, it has been shown in \cite{DPD-NS3D}
that it is possible to construct a Markov transition semigroup
which is the limit of Galerkin approximations for \eqref{e1.1}.
Moreover, this semigroup is associated to weak solutions and thanks
to the nondegeneracy of the noise, it  is strongly mixing. It
follows from \cite{ODASSO6} that it is exponentially mixing.

Also, in \cite{flandoli-cetraro}, \cite{romito-flandoli-selection}, another idea of
the book \cite{SV} has been used. It has been
shown that it is possible to construct a Markov selection from weak solutions to \eqref{e1.1}.
Moreover, using the ideas of \cite{DPD-NS3D}, it is proved that if
the noise is sufficiently non degenerate
every Markov selection of weak solutions defines a Strong Feller transition semigroup.

In this article, our aim is to  improve the result of
\cite{DPD-NS3D} and to extend it to the case of a state dependent
noise. Moreover, we considerably simplify the construction.
Indeed, in \cite{DPD-NS3D}, numerous technical a priori estimates
are necessary for a suitable approximation of the Kolmogorov
equation introduced thanks to Galerkin approximation. We use the
same idea here but show that it is sufficient to estimate the
first derivate of their solutions as well as their modulus of
continuity in time. This allows to prove compactness of the
approximated solution of the Kolmogorov equations. We then notice
that in fact this implies the convergence in law of a subsequence of the
family of solutions. The limit is clearly a Markov family of
solutions. Contrary to \cite{flandoli-cetraro},
\cite{romito-flandoli-selection}, this Markov family of solution
is obtained in a constructive way.

We also prove irreducibility of the transition semigroup thanks to
classical arguments taken from the method to prove support
theorems. Ergodicity follows.

\section{Preliminaries}
    We set
$$
H=\{x\in (L^2(\mathcal O))^3:\;\mbox{\rm div}\; x=0\;\mbox{\rm
in}\;\mathcal O,\;x\cdot n=0\;
\mbox{\rm on}\;\partial \mathcal O\},
$$
where $n$ is  the outward normal to $\partial \mathcal O,$
and $V= (H^1_0(\mathcal O))^3\cap H.$
The norm and  inner product in $H$   will be denoted by $|\cdot|$   and
 $(\cdot,\cdot)$ respectively. We denote by $\Lr(H)$ (resp. $\Lr_2(H)$) 
 the space of linear (resp. Hilbert-Schmidt) operators on $H$ with 
 norm $|\cdot|_{\Lr}$ (resp.  $|\cdot|_{\Lr_2}$).
Moreover $W$ is a cylindrical Wiener process  on $H$ and, for any
$x\in H$, the  operator $\Phi(x) \in \Lr(H)$ is Hilbert-Schmidt
 and such that Ker $\phi(x)=\{0\}$.

  \vspace{2mm}
 \noindent We also denote by $A$ the Stokes operator in $H$:
    $$A= P\Delta,\quad D(A)=(H^2(\mathcal O))^3\cap (H^1_0(\mathcal
    O))^3
\cap H,$$
where   $P$ is the orthogonal projection of $(L^2(\mathcal O))^3$ onto $H$
and by $b$  the operator
$$
b(x,y) =-P((x\cdot\nabla)y),\quad
b(x)=b(x,x),\quad x,y\in V.
$$
Now we can write problem \eqref{e1.1} in the form
\begin{equation}
\label{e1.2}
\left\{\begin{array}{l}
dX(t,x)= (AX(t,x)+b(X(t,x))+f)dt+\Phi(X(t,x))\;dW(t),\\
\\
X(0,x)=x.
\end{array}\right.
\end{equation}
We assume that the forcing term $f$ is in $V$.

Classically, we use fractional powers of the operator $A$ as
well as their domains $D((-A)^\alpha)$ for $\alpha \in \R$.
Recall that, thanks to the regularity theory of the Stokes operator,
$D((-A)^\alpha)$ is a closed subspace of the Sobolev space
$(H^{2\alpha}({\cal O}))^{3}$ and $|\cdot |_{D((-A)^\alpha)} =
|(-A^\alpha)\cdot|$ is  equivalent to the usual
$(H^{2\alpha}({\cal O}))^{3}$ norm.

Let $0<\alpha<\beta<\gamma.$ Then the following interpolatory estimate is well known
\begin{equation}
    \label{einter}
|(-A)^\beta x|\le c |(-A)^\alpha x|^{\frac{\gamma-\beta}{\gamma-\alpha}}\;|(-A)^\gamma x|^{\frac{\beta-\alpha}
{\gamma-\alpha}}, \quad x\in D((-A)^\gamma).
\end{equation}
Moreover
$
(b(x,y),y)=0,
$
whenever the left hand side makes sense. We shall use the following estimates on  the bilinear operator $b(x,y)$ (see
\cite{constantin-foias},  \cite{DPD-NS3D}, \cite{Temam}):
\begin{equation}
\label{interpolation} \left|(-A)^{1/2}b(x,y)\right|\le   c |Ax|\;
|Ay|.
\end{equation}

The following  result can be proved by classical argument (see
\cite{flandoli-cetraro}).
\begin{Proposition}
\label{p1.1}
For any $x\in H$, there exists a martingale solution of equation
\eqref{e1.2} with trajectories in   $C([0,T];D((-A)^{-\alpha}))$
and $L^\infty(0,T;H)\cap L^2(0,T;D((-A)^{1/2})$
for any $\alpha>0$ and $T> 0$.
\end{Proposition}
Let $E$ be  any Banach space and  $\varphi:D(A)\to E.$ For any
$x,h\in D(A)$ we set
$$
D\varphi(x)\cdot h=\lim_{s\to 0}\frac 1s\;(\varphi(x+sh)-\varphi(x)),
$$
provided the limit exists.  The limit is intended in $E.$ The  functional space
$C_b(D(A);E)$  is the space of all continuous and bounded
mappings from $D(A)$ (endowed with the graph norm)
into $E$ and, for any $k\in \N$, $C_k(D(A);E)$  is the space of all
continuous  mappings from $D(A)$ into $E$
such that
$$
\|\varphi\|_{k}:=\sup_{x\in D(A)}\frac{|\varphi(x)|}{(|Ax|+1)^k}<+\infty.
$$
We denote by $B_b(D(A);\R)$ the space of all Borel bounded
mappings from $D(A)$ into $\R$.
It is also convenient to define the space
$$
{\mathcal E}=\left\{\varphi\in
C_b(D(A);\R)\,\left|\,
\sup_{(x_1,x_2)}\frac{|\varphi(x_2)-\varphi(x_1)|}{|A(x_2-x_1)|(1+|Ax_1|^2+|Ax_2|^2)}
<\infty\right.\right\},
$$
with the norm
$$
\|\varphi\|_{{\mathcal E}} =\|\varphi\|_0+\sup_{(x_1,x_2)}\frac{|\varphi(x_2)-\varphi(x_1)|}{|A(x_2-x_1)|(1+|Ax_1|^2+|Ax_2|^2)}
$$
We introduce the usual Galerkin approximations of equations
\eqref{e1.2}. For $m\in \N$, we define the projector $P_m$ onto
the first $m$ eigenvectors of $A$ and set $b_m(x)=P_mb(P_m x)$,
$\Phi_m(x)=P_m \Phi(x)$, for $x\in H$. Then, we write the
following approximations
\begin{equation}
\label{e1.4}
\left\{\begin{array}{l}
dX_m(t,x)= (AX_m(t,x)+b_m(X_m(t,x)))dt+\Phi_m(X_m(t,x))\;dW(t),\\
\\
X_m(0,x)=P_mx=x_m,
\end{array}\right.
\end{equation}
and
\begin{equation}
\label{e1.5}
\left\{\begin{array}{l}
\ds{\frac{du_m}{dt}\;=\frac 12\;\mbox{\rm Tr}\;[\psi_m(x)D^2u_m]+(Ax+b_m(x),Du_m)},\\
\\
u_m(0)=\varphi,
\end{array}\right.
\end{equation}
where $\psi_m=\Phi_m\Phi_m^*$.
Equation \eqref{e1.5} has a unique solution given by
\begin{equation}
\label{e1.6}
u_m(t,x)=P^m_t\varphi(x)=\E[\varphi(X_m(t,x))],\;\textrm{ for
}x\in P_m H.
\end{equation}
We extend $u_m(t,x)$ to $H$ by setting $u_m(t,x)=u_m(t,P_m x)$.
We assume that $\Phi$ is a $C^1$ function of $x$ with
values in $\Lr(H)$. Then, if $\varphi$ is a $C^1$ function, $u_m$ is
differentiable and its differential can be expressed in terms of
 $\eta_m^h=\eta_m^h(t,x)$  the solution of
\begin{equation}
\label{e1.8}
\left\{\begin{array}{l}
\ds{d\eta_m^h=
\left(A\eta_m^h+b'_m(X_m)\cdot\eta_m^h\right)dt
+ \Phi_m'(X_m)\cdot \eta_m^h dW(t)},\\
\\
\eta_m^h(0,x)=P_m h,
\end{array}\right.
\end{equation}
with  $b'_m(X_m)\cdot\eta_m^h=b_m(X_m,\eta_m^h)
+b_m(\eta_m^h,X_m).$ Moreover, since we assume that the noise is
non degenerate, the differential of $u_m$ exists even when
$\varphi$ is only continuous thanks to the Bismut--Elworthy--Li
formula, see \cite{B} and \cite{E}. Unfortunately, it is
impossible to  get any estimate on the differential of $u_m$ by
these ideas. Indeed, we are not able to prove an estimate of
$\eta_m^h(t,x)$ uniform in $m$. The idea is to introduce the
auxiliary Kolmogorov equation
\begin{equation}
\label{e1.10}
\left\{\begin{array}{l}
\ds{\frac{dv_m}{dt}\;=\frac 12\;\mbox{\rm
Tr}\;[\psi_m(x)D^2v_m]+(Ax+b_m(x),Dv_m)-K|Ax|^2v_m},\\
\\
v_m(0)=\varphi,
\end{array}\right.
\end{equation}
where $K>0$ is fixed, which contains a ``very negative''  potential term.
It has a unique solution given by the Feynman-Kac formula
\begin{equation}
\label{e1.11}
v_m(t,x):=S^m_t\varphi(x)=\E\left[e^{-K\int_0^t|AX_m(s,x)|^2ds}\varphi(X_m(t,x))\right].
\end{equation}
Clearly, the function $u_m$ can be expressed in terms of the function
$v_m$ by the variation
of constants formula:
\begin{equation}
\label{e1.12}
u_m(t,\cdot)=S^m_t\varphi+K\int_0^tS^m_{t-s}(|A\cdot|^2u_m(s,\cdot))ds.
\end{equation}
Since the noise is non degenerate, from \cite{DPZ3} we know that for any
$\varphi\in C_b(H),$
$S^m_t\varphi$ is differentiable in any
direction
$h\in H$ and we have
\begin{equation}
\label{e1.14}
\begin{array}{l}
DS^m_t\varphi(x)\cdot h\\
\\
\ds{=\frac{1}{t}\;\E\left[e^{-K\int_{0}^{t}|AX_m(s,x)|^2ds}\varphi(X_m(t,x))
 \int_{0}^{t}(\Phi^{-1}(X_m(s,x))\eta_m^h(s,x),dW(s))
\right]}\\
\\
\ds{+2K\E\left[e^{-K\int_{0}^{t}|AX_m(s,x)|^2ds}\varphi(X_m(t,x))
  \int_{0}^{t}\left( 1-\frac{s}{t} \right)(AX_m(s,x),A\eta_m^{h}(s,x))ds \right].}\\
\\
\end{array}
\end{equation}
We are going to prove estimates for the  derivatives of
$u_m(t,\cdot)$ through
corresponding estimates for $v_m(t,\cdot)$. We will see that this is possible provided
$K$ is chosen large enough. This implies some compactness
on the sequence $(u_m)_m$.

The main assumption in our estimates below is that the covariance operator
is at the same
time sufficiently smooth and non degenerate. We assume throughout the paper that there
exists  constants $M_1\ge 0$, $r\in (1,3/2)$ and $g>0$ such that
\begin{equation}
\label{e1.18}
\mbox{\rm Tr}\;[(-A)^{1+g}\Phi(x)\Phi^*(x)]\le M_1,\mbox{ for }\; x \in H
\end{equation}
   and
\begin{equation}
\label{e1.19}
|\Phi^{-1}(x)h|\le M_1 |(-A)^{r}h|,\quad\mbox{for } x\in H, \,
h\in D((-A)^{r}).
\end{equation}
   Note that these two
conditions are compatible, for instance if we take $\Phi(x)$ which
is bounded and invertible from $H$ onto $D((-A)^{\alpha})$,
$\alpha>0$ then conditions \eqref{e1.18} and \eqref{e1.19} are
satisfied provided
$$
\alpha \in (5/4,3/2),
$$
and the norm of $\Phi(x)$ and its inverse are uniformly bounded.

Since we also work with the differential of the solution with respect to the initial data, we also
need some assumption on the derivative of $\Phi$. We assume that there exists $\delta<3/2$ such
that:
\begin{equation}
\label{e1.20}
\mbox{\rm Tr}\;[(-A)^{2}(\Phi'(x)\cdot h)(\Phi'(x)\cdot h)^*]\le M_1|(-A)^\delta h|^2,\mbox{ for }\; x \in H, \;
h\in D((-A)^\delta).
\end{equation}

\begin{Example}
Let $c>0$, $\alpha \in (5/4,3/2)$ and $\kappa$ a bounded  mapping
$H\to {\mathcal L}_2 (H;D(A))$ with bounded continuous
derivative. We set
$$
\Phi(x)=(-A)^{\alpha}+c\kappa(x).
$$
If Ker $\Phi(x)= \{0\}$ for any $x\in H$, then $\Phi$ verifies the previous assumptions.
This is the case if $c$ is
small.
A first  example of such $\kappa$ is the following. Let $(\kappa_n)_n$
be a family of functions in $C^1_b(H)$ with norm bounded by one  and
$(\lambda_n)_n\in l^2(\N)$. We denote by $(\mu_n,e_n)$ the family
of eigenvalues and eigenvectors of $(-A)$ . Then,
we set
$$
\kappa(x)\cdot h=\sum_n\mu_n^{-\frac{3}{2}} \lambda_n\kappa_n(x)
h_ne_n,\quad \textrm{ where } h=\sum_n h_ne_n.
$$
Another example is given by
$$
(\kappa(x)\cdot h)(\xi)=\int_{O\times O}\mathcal
V(\xi,\xi',x(\xi'))h(\xi')\,d\xi',
$$
where  $\mathcal V$ is a $C^\infty$ map $O\times O\times\R\to \R$
with compact support and free divergence with respect to the first
variable (i.e. div $\mathcal V(\cdot,\xi,r)=0$).
\end{Example}
Before stating our main result, we recall the following definition
\begin{Definition}\label{Def}
Let $(\Omega_x,\mathcal F_x,\P_x)_{x\in D(A)}$ be a family of
probability spaces and $(X(\cdot,x))_{x\in D(A)}$ be a family of
random processes on $(\Omega_x,\mathcal F_x,\P_x)_{x\in D(A)}$.
We denote by
$(\mathcal F_x^t)_t$ the filtration generated by $X(\cdot,x)$ and  by
$\mathcal P_x$ the law of $X(t,x)$ under
$\P_x$.
The family $(\Omega_x,\mathcal
F_x,\P_x,X(\cdot,x))_{x\in D(A)}$ is a Markov family if
the following conditions hold:
\begin{enumerate}
\item[(i)] For any $x\in D(A)$, $ t\ge 0$, we have
$$
\P_x\left(X(t,x)\in D(A)\right)=1,
$$

 \item[(ii)] the
map $x\to \mathcal P_x$ is measurable and for any $x\in D(A)$,
$t_0,\dots,t_n\ge 0$,
$A_0,\dots ,A_n$ borelian subsets of $D(A)$, we have
$$
\P_x\left(X(t+.)\in  {\mathcal A}\,\left |\, \mathcal
F_x^t\right.\right)=\mathcal P_{X(t,x)}(\mathcal A),
$$
where ${\mathcal A}=\{y\in (H)^{\R_+}\,|\,  y(t_0)\in  A_0,\dots, y(t_n)\in A_n\}$.
\end{enumerate}
The Markov transition semigroup  $(P_t)_{t\ge 0}$ associated to
the family is then defined by
$$
P_t \varphi(x) = \E_x(\varphi(X(t,x))), \quad x\in D(A),\; t\ge 0.
$$
for $\varphi \in B_b(D(A);\R)$.
\end{Definition}
The main result of the paper is the following. The proof is given in the case $f=0$, the generalization
to a general $f\in V$ is  easy.
\begin{Theorem}
\label{t7.1} There exists a Markov family of martingale solutions
$(\Omega_x,\mathcal F_x,\P_x,X(\cdot,x))_{x\in D(A)}$ of the
stochastic Navier-Stokes equations \eqref{e1.2}.
Furthermore, the transition semigroup $(P_t)_{t\ge 0}$ is
stochastically continuous.
\end{Theorem}
We also study ergodic properties and prove the following result.
\begin{Theorem}
\label{t7.2}
 There exists a Markov process $X(\cdot,\nu)$ on a probability 
 space  $(\Omega_\nu,\mathcal F_\nu,\P_\nu)$ which is a martingale 
 stationary solution
 of the
stochastic Navier-Stokes equations \eqref{e1.2}. The law
$\nu$ of $X(t,\nu)$ is the unique invariant measure on $D(A)$ of
the transition semigroup $(P_t)_{t\ge 0}$. Moreover
\begin{enumerate}
\item[(i)] the invariant measure $\nu$ is ergodic,

\item[(ii)]  the law  $\mathcal P_\nu$ of
$X(\cdot,\nu)$
 is given by
 $$
\mathcal P_\nu(\mathcal A)=\int_{D(A)}\mathcal P_x(\mathcal A)\,\nu(dx),
 $$
for  ${\mathcal A}=\{y\in (H)^{\R_+}\,|\,  y(t_0)\in  A_0,\dots, y(t_n)\in A_n\}$
with $t_0,\dots,t_n\ge 0$ and
$A_0,\dots ,A_n$ borelian subsets of $D(A)$,
\item[(iii)] the transition semigroup $(P_t)_{t\ge 0}$ is strong
Feller, irreductible, and therefore strongly mixing.
\end{enumerate}
\end{Theorem}
\begin{Remark} In fact, we prove that
there exists a subsequence $(m_k)_k$ such that, for any $x\in D(A)$,
$X_{m_k}(\cdot,x)\to X(\cdot,x)$ in law.
Moreover $X_{m_k}(\cdot,\nu_{m_k})$ the unique stationary
solution of \eqref{e1.4} converges to $X(\cdot,\nu)$. Thus, the
solutions $(X(\cdot,x))_x$ and $X(\cdot,\nu)$ are obtained in a
constructive way.
\end{Remark}

\section{A priori estimates}

For any predictible process $X$ with values in $H$, we set:
\begin{equation}
    \label{e1.20bis}
Z_X(t) =\int_0^t e^{(t-s)A}\Phi(X(s))\;dW(s).
\end{equation}
We have the following estimates on $Z_X$ which will be useful in the sequel.
\begin{Proposition}\label{P2.2}
For  any $T\ge 0$, $\vep<g/2$ and any $m\ge 1$, there exists a constant $c(\vep,m,T)$ such that,
for any predictible process $X$ with values in $H$, $Z_X$ has continuous paths with values
in $D((-A)^{1+\vep})$ and
\begin{equation}
\label{2.15}
\E(\sup_{t\in [0,T]}|(-A)^{1+\vep}Z_X(t)|^{2m})\le c(\vep,m,T).
\end{equation}
Moreover, for any $\beta < \min\{g/2 -\vep,1/2\}$, there exists a
constant $c(\vep,\beta,m,T)$ such that for $t_1,t_2\in [0,T]$,
\begin{equation}
\label{2.16}
\E(|(-A)^{1+\vep}(Z_X(t_1)-Z_X(t_2))|^{2m})\le c(\vep,\beta,m,T)|t_1-t_2|^{2\beta m}.
\end{equation}

\end{Proposition}
{\bf Proof:} The proof uses the factorization method (see \cite[Section 5.3]{DPZ1}). We write
$$
Z_X(t)=\int_0^t(t-s)^{\alpha-1}e^{A(t-s)}Y(s)ds
$$
with $\alpha$ to be chosen below and
$$
Y(s)=\frac{\sin \pi\alpha}\alpha\int_0^s e^{A(t-s)}(s-\sigma)^{-\alpha}\Phi(X(\sigma))dW(\sigma).
$$
Using Burkholder-Davies-Gundy inequality, we deduce for $m\in \N$:
$$
\begin{array}{l}
\ds{\E\left(|(-A)^{1+\vep}Y(s)|^{2m}\right)}\\
\ds{\le c \E\left(\left(\int_0^s  |(-A)^{1+\vep} e^{A(s-\sigma)}
(s-\sigma)^{-\alpha}\Phi(X(\sigma))|^2_{\Lr_2}d\sigma\right)^m\right)}\\
\ds{\le c \E\left( \left( \int_0^s (s-\sigma)^{-2\alpha} |(-A)^{1/2+g/2}\Phi(X(\sigma))|_{\Lr_2}^2 |(-A)^{1/2+\vep-g/2}e^{A(s-\sigma)}|_{\Lr}^2 d\sigma\right)^m\right)}\\
\ds{\le c M_1 \left(\int_0^s
(s-\sigma)^{-(1+2\vep-g)^+-2\alpha}d\sigma\right)^m,}
\end{array}
$$
thanks to well-known smoothing properties of the semigroup
$(e^{At})_{t\ge 0}$ and assumption \eqref{e1.18}. This is a finite
quantity provided $\alpha<\min\{g/2-\vep,1/2\}$. Moreover, it is a
bounded function of $s\in [0,T]$.
 It follows easily that $Y\in L^{2m}(\Omega\times [0,T]; D((-A)^{1+\vep}))$.
 Proposition \ref{P2.2} follows now from Proposition
A.1.1 of \cite{DPZ2}.  \hfill $\Box$

The proof of the following estimate  is  the same as the proof of Lemma 3.1
in \cite{DPD-NS3D}.
\begin{Lemma}
\label{l2.1}
There exists $c>0$  such
that, for any $m\in \N$, $t\in [0,T]$ and $x\in D(A)$,
\begin{equation}
\label{e2.2c} e^{-c\int_0^t |AX_m(s,x)|^2 ds}|AX_m(t,x)|^2 \le
2|Ax|^2+c \sup_{s\in [0,T]} |AZ_{X_m}(s)|^2.
\end{equation}
\end{Lemma}
\begin{Lemma}
\label{l2.2}
For any $\gamma \in (\delta-1/2, 1]$, there exists $c_\gamma>0$  such
that for any $m\in \N$, $t\in [0,T]$ and any
$x,h\in D(A)$ we have
\begin{equation}
\label{e7.1}
\begin{array}{l}
\ds{\E\left(e^{-c_\gamma\int_0^t |AX_m(s,x)|^2 ds}|(-A)^\gamma\eta_m^h(t,x)|^2
+\int_0^t e^{-c_\gamma\int_0^s |AX_m(\tau,x)|^2 d\tau}|(-A)^{\gamma+1/2}\eta_m^h(s,x)|^2ds\right)}\\
\ds{\le e^{c_\gamma t}  |(-A)^\gamma h|^2,}
\end{array}
\end{equation}
\end{Lemma}
{\bf Proof}. We use Ito formula to obtain
$$
\begin{array}{l}
\ds{d\bigg[e^{-c_\gamma \int_0^t |AX(s)|^2 ds}|(-A)^\gamma\eta(t)|^2
+2\int_0^t e^{-c_\gamma\int_0^s |AX(\tau)|^2 d\tau}|(-A)^{\gamma+1/2}\eta(s)|^2ds\bigg]}\\
\ds{= e^{-c_\gamma \int_0^t |AX(s)|^2 ds}\bigg[
\mbox{\rm Tr}\left[ (-A)^{2\gamma}\left(\Phi'(X(t))\cdot
\eta(t)\right) \left(\Phi'(X(t))\cdot \eta(t)\right)^*\right]dt}\\
\ds{+2\left((-A)^{2\gamma}\eta(t),\left(\Phi'(X(t))\cdot
\eta(t)\right)dW(t)\right)+
2\left(b'(X(t))\cdot\eta(t),(-A)^{2\gamma}\eta(t)\right)dt}\\
\ds{-c_\gamma |A X(t)|^2|(-A)^{\gamma}\eta(t)|^2dt
\bigg].}
\end{array}
$$
We have written for simplicity $\eta_m^h(t)=\eta(t)$ and $X(t)=X_m(t,x)$.
Using \eqref{interpolation}, Poincar\'e inequality and \eqref{einter}, we obtain
$$
\begin{array}{ll}
\ds{\left(b'(X(t)),\eta(t),(-A)^{2\gamma}\eta(t)\right)}
&\ds{\le c\;|A X(t)| \; |A \eta(t)|\; |(-A)^{2\gamma-1/2}\eta(t)|}\\
&\ds{\le c \;  |A X(t)| \; |(-A)^\gamma \eta(t)|\; |(-A)^{\gamma+1/2}\eta(t)|}\\
&\ds{\le c \;  |A X(t)|^2|(-A)^\gamma \eta(t)|^2+\frac 14 |(-A)^{\gamma+1/2}\eta(t)|^2}.
\end{array}
$$
Moreover, since $\delta<\gamma +1/2$ in \eqref{e1.20}, we have
$$
\begin{array}{rl}
\lefteqn{\ds{\mbox{\rm Tr}\left[
(-A)^{2\gamma}\left(\Phi'(X(t))\cdot \eta(t)\right)
\left(\Phi'(X(t))\cdot \eta(t)\right)^*\right] }}\\
&\ds{\le c\; \mbox{\rm Tr}\left[ (-A)^{2}\left(\Phi'(X(t))\cdot
\eta(t)\right)
\left(\Phi'(X(t))\cdot \eta(t)\right)^*\right] }\\
&\ds{\le c\; |(-A)^\delta \eta(t)|^2}\\
&\ds{\le c \; |(-A)^\gamma \eta(t)|^2+\frac
12|(-A)^{\gamma+1/2}\eta(t)|^2}.
\end{array}
$$
It follows that, for $c_\gamma$ sufficiently large,
$$
\begin{array}{l}
\ds{d\bigg[e^{-c_\gamma \int_0^t |AX(s)|^2 ds}|(-A)^\gamma\eta(t)|^2
+\int_0^t e^{-c_\gamma\int_0^s |AX(\tau)|^2 d\tau}|(-A)^{\gamma+1/2}\eta(s)|^2ds\bigg]}\\
\le \ds{ 2e^{-c_\gamma \int_0^t |AX(s)|^2
ds}\bigg[\left((-A)^{2\gamma}\eta(t),\left(\Phi'(X(t))\cdot
\eta(t)\right)dW(t)\right) +c \; |(-A)^\gamma \eta(t)|^2dt\bigg]}.
\end{array}
$$
We deduce the result by taking the expectation and integrating. \hfill $\Box$

We now get bounds on the Feynman-Kac semigroup $S_t^m$.
\begin{Lemma}
\label{l3.1} For any $1\ge \gamma>\max \{\delta -1/2, r-1/2\}$,
where $r\in (1,3/2)$ is defined in \eqref{e1.19} and $k\in \N$, if
$K $ is sufficiently large there exists $c(\gamma,k)>0$ such that
for any $\varphi \in C_k(D(A);\R)$
$$
\|(-A)^{-\gamma} DS_t^m\varphi\|_{k} \le c (\gamma)
(t^{-1/2-(r-\gamma)}+1)\|\varphi\|_{k},\quad t>0,
$$
for all $m\in \N.$
\end{Lemma}
{\bf Proof}.
Let $h \in H.$ We write
\eqref{e1.14} as  $DS^m_t\varphi(x)\cdot h=I_1+I_2$ and estimate separately the two terms.
We again write for simplicity $\eta_m^h(t)=\eta(t)$ and $X(t)=X_m(t,x)$.
Concerning  $I_1$ we have, using the H\"older inequality,
$$
    \begin{array}{l}
\ds{I_1 \le \frac{1}{t}\;
\|\varphi\|_{k}\;\E\left[e^{-K\int_{0}^{t}|AX(s)|^2ds}(1+|AX(t)|)^k
\left(\int_{0}^{t}(\Phi^{-1}(X(s))\eta(s),dW(s))\right)\right]}\\
\\
\ds{\le
   \frac{1}{t}\;
\|\varphi\|_{k}\;\left[\E\left(e^{-K\int_{0}^{t}|AX(s)|^2ds}(1+|AX(t)|)^{2k}
\right)\right]^{1/2}}\\
\\
\ds{ \times\left[\E\left(e^{-K\int_{0}^{t}|AX(s)|^2ds}
\left(\int_{0}^{t}(\Phi^{-1}(X(s))\eta(s),dW(s))\right)^2\right)\right]^{1/2}.}
    \end{array}
$$
Choosing $K$ sufficiently large, the first factor is easily majorized by $c (1+|Ax|)^{2k}$ thanks to Lemma \ref{l2.1} and Proposition \ref{P2.2}.
To estimate the second factor
we proceed as in \cite{DPD-NS3D} and set
$$
\xi(t)=e^{-\frac{K}{2}\int_{0}^{t}|AX(s)|^2ds}
\int_{0}^{t}(\Phi^{-1}(X(s))\eta(s),dW(s))
$$
and use Ito formula to compute $\E(\xi^2(t))$. We obtain
$$
\E(\xi^2(t))\le \E\left[\int_0^te^{-K\int_{0}^{s}|A
X_m(\tau,x)|^2d\tau}|\Phi^{-1}(X(s))\eta(s)|^2ds \right].
$$
Recalling   assumption \eqref{e1.19} and the interpolatory estimate   \eqref{einter}
 we find
$$
|\Phi^{-1}(X(s))x|^2\le M_1 |(-A)^rx|^2\le
c|(-A)^\gamma x|^{2(1-2(r-\gamma))}\;|(-A)^{\gamma +1/2}x|^{4(r-\gamma)}.
$$
Consequently, by H\"older inequality and    Lemma \ref{l2.2},
we get
$$
\E(\xi(t)^2)\le
 c t^{1-2(r-\gamma)}|(-A)^\gamma h|^2,
$$
provided $K$ is sufficiently large.
Thus
$$
I_1\le ct^{-1/2-(r-\gamma)}|(-A)^\gamma h|\;(|Ax|+1)^k.
$$
Finally, since $\gamma \ge 1/2$, the following estimate easily follows from
H\"older inequality and  Lemmas
\ref{l2.1}, \ref{l2.2}:
$$
I_2\le c\|\varphi\|_{k}\;(1+|Ax|)^k|(-A)^\gamma h|.
$$
Consequently, if
$K$ is sufficiently large,  we find
$$
|DS_t^m\varphi(x)\cdot h|\le c\|\varphi\|_{k}(1+|Ax|)^k
(1+t^{-1/2-(r-\gamma)})|(-A)^\gamma h|.
$$
The conclusion follows.
\hfill $\Box$

We are  now ready   to get uniform
estimates on the approximated solutions to the Kolmogorov equation.
\begin{Proposition}
\label{p4.1} If $\varphi\in C_b(D(A);\R),$ then $u_m(t)\in
C_b(D(A);\R)$ and,  for any $1\ge \gamma> \max \{\delta -1/2,
r-1/2\}$, $(-A)^{-\gamma}Du_m \in C_2(D(A);\R) $ for all $t>0$,
$m\in \N$. Moreover, we have
$$
\|u_m(t)\|_0\le \|\varphi \|_0
$$
and
$$
\|(-A)^{-\gamma}Du_m(t)\|_{2}\le c(\gamma)
(1+t^{-1/2-(r-\gamma)})\|\varphi\|_{0},
\quad t>0.
$$
\end{Proposition}
{\bf Proof}.
The first estimate follows from \eqref{e1.6}. By \eqref{e1.12} and Lemma \ref{l3.1}, it follows
$$
\begin{array}{ll}
\ds{\|(-A)^{-\gamma}Du_m(t)\|_{2}}&\ds{\le c
(1+t^{-1/2-(r-\gamma)})\|\varphi\|_{2}}\\
&\ds{+\int_0^t c(1+(t-s)^{-1/2-(r-\gamma)})\||Ax|^2u_m(s)\|_{2}\;ds.}
\end{array}
$$
Clearly $\|\varphi\|_{2}\le \|\varphi\|_0$ and
$
\||Ax|^2u_m(s)\|_{2}\le \|u_m(s)\|_{0}\le \|\varphi \|_{0}.
$
The result follows. \hfill $\Box$

\begin{Proposition}
\label{pc}
Let  $\varphi\in   {\mathcal E}$. Then for any $\beta <\min\{g/2,1/2\}$, there exists $c(\beta)$ such that for
any
 $t_1,t_2>0$, $m\in \N$ and $x\in D(A)$ we have
$$
|u_m(t_1,x)-u_m(t_2,x)|\le
c\|\varphi\|_{\mathcal E} \left(|Ax|+1\right)^6\left(|t_1-t_2|^{\beta}
+
|A(e^{t_1A}-e^{t_2A})x|\right).
$$
\end{Proposition}
{\bf Proof}. By \eqref{e1.12}, we have
 for $t_1<t_2$
$$
\begin{array}{ll}
\ds{u_m(t_1,x)-u_m(t_2,x)}&\ds{=\left(S^m_{t_1}-S^m_{t_2}\right)\varphi (x)
+K\int_0^{t_1}\left(S^m_{t_1-s}-S^m_{t_2-s}\right)
\left(|Ax|^2 u_m(s)\right) (x)ds
}\\
\\
&\ds{+K\int_{t_1}^{t_2}S^m_{t_2-s}\left(|Ax|^2 u_m(s)\right) (x)ds
}\\
\\
&=T_1+T_2+T_3.
\end{array}
$$
For the first term  we use the decomposition, with $X(t)=X_m(t,x)$,
$$
\begin{array}{ll}
\ds{|T_1| } &\ds{=\bigg| \E\left( \left(e^{-K\int_0^{t_1}|AX(s)|^2ds}-e^{-K\int_0^{t_2}|AX(s)|^2ds}
\right)\varphi(X(t_1))\right)}\\
\\
&\ds{ + \E\left(e^{-K\int_0^{t_2}|AX(s)|^2ds}
\left(\varphi(X(t_1))-\varphi(X(t_2)\right)\right) \bigg|}\\
\\
&\ds{\le K \|\varphi\|_{\mathcal E} \E\left( \int_{t_1}^{t_2} |AX(s)|^2 e^{-K\int_0^{s}|AX(\tau)|^2d\tau} ds\right)}\\
\\
&\ds{+ \|\varphi\|_{\mathcal E} \E\left(e^{-K\int_0^{t_2}|AX(s)|^2ds}(1+ \sup_{t\in [0,T]} |AX(t)|^2)
|A\left(X(t_1)-X(t_2)\right)|\right) }.
\end{array}
$$
Thanks to Lemma \ref{l2.1}, we majorize the first term by
$c \|\varphi\|_0 (|Ax|^2+1)|t_1-t_2|$. For the second term, we write
$$
\begin{array}{ll}
X(t_1)-X(t_2)&\ds{=(e^{t_1A}-e^{t_2A})x+Z_{X}(t_1)-Z_{X}(t_2)}\\
&\ds{+\int_0^{t_1}e^{A(t_1-s)}b(X(s))ds
-\int_0^{t_2}e^{A(t_2-s)}b(X(s))ds.}
\end{array}
$$
By  \eqref{interpolation}, and classical property of $(e^{At})_{t>0}$, for any $\lambda \in (0,1/2)$ we have the following estimate:
$$\begin{array}{l}
\ds{\E\left(e^{-K\int_0^{t_2}|AX(s)|^2ds}
\left|A\left(\int_0^{t_1}e^{A(t_1-s)}b(X(s))ds
-\int_0^{t_2}e^{A(t_2-s)}b(X(s))ds  \right)\right|\right)}\\
\ds{\le \E\bigg[e^{-K\int_0^{t_2}|AX(s)|^2ds}\bigg(\int_{t_1}^{t_2}\left|(-A)^{1/2}e^{A(t_2-s)}\right|_{\Lr}
\left| (-A)^{1/2} b(X(s))\right| ds}\\
\ds{\quad\quad\quad\quad+\int_0^{t_1}
\left|(-A)^{1/2}\left(e^{A(t_1-s)} -
e^{A(t_2-s)}\right)\right|_{\Lr}
 \left| (-A)^{1/2}b(X(s))\right|ds\bigg)\bigg]}\\
\ds{\le c_\lambda
\E\left(e^{-K\int_0^{t_2}|AX(s)|^2ds}\sup_{s\in[0,t_2]} \left|
AX_m(s,x)\right|^2\right) } |t_1-t_2|^{1/2-\lambda}
\end{array}
$$
Therefore, with  Proposition \ref{P2.2} and Lemma \ref{l2.1}, we
obtain for any $\beta<\min\{ g/2, 1/2\}$
$$
|T_1|\le c\|\varphi\|_{\mathcal E}(|Ax|^4+1)\left(|t_1-t_2|^{\beta}
+\left|A(e^{t_1A}-e^{t_2A})x\right|\right).
$$
Similarly, we have
$$
\begin{array}{ll}
\left(S^m_{t_1-s}-S^m_{t_2-s}\right)
\left(|Ax|^2 u_m(s)\right)& \le c(\|u_m\|_0+\|(-A)^{-1}Du_m(s)\|_{2})(|Ax|^6+1)\\
\\
&\times \left(|t_1-t_2|^{\beta}
+\left|A(e^{t_1A}-e^{t_2A})x\right|\right).
\end{array}
$$
Thus, thanks to Proposition \ref{p4.1},
$$
|T_2|\le c \|\varphi\|_{\mathcal E}(|Ax|^6+1)\left(|t_1-t_2|^{\beta}
+\left|A(e^{t_1A}-e^{t_2A})x\right|\right).
$$
Finally, the last term $T_3$ is easy to treat and majorized by $c \|\varphi\|_0 (|Ax|^2+1)|t_1-t_2|$.
Gathering these estimates yields the result. \hfill $\Box$

\section{Proof of Theorem \ref{t7.1}}

For $\varphiÊ\in {\mathcal E}$, let $(u_m)_{m\in \N}$ be the sequence of solutions of
the approximated Kolmogorov equations. Thanks to the a priori estimates derived in the
previous section, we are now able to show that $(u_m)_{m\in \N}$ has a convergent
subsequence.

We set $K_R=\{x\in D(A):\;|Ax|\le R\}$. For $\gamma<1$, it   is a compact
subset of $D((-A)^\gamma).$
\begin{Lemma}
\label{l5.1}
Assume that $\varphi \in {\mathcal E}$, then there exists a subsequence
$(u_{m_k})_{k\in \N}$ of $(u_{m})$ and a function $u$  bounded on $[0,T]\times D(A),$   such that
\begin{enumerate}
\item[(i)]
$u\in C_b((0,T]\times D(A))$ and for any $\delta>0$, $R>0$
$$
\lim_{k\to \infty}u_{m_k}(t,x)=u(t,x)\quad\mbox{\it uniformly
on}\;[\delta,T]\times K_R.
$$
\item[(ii)] For any $x\in
D(A)$, $u(\cdot,x)$ is continuous on $[0,T]$.
\item[(iii)] For any
$1\ge \gamma> \max \{\delta -1/2, r-1/2\}$, $\delta >0$, $R\ge 0$,
$\beta<\min\{g/2,1/2\}$, there exists
 $c(\gamma,\beta,\delta, R, T, \varphi)$
such that for $x,y\in K_R$, $t,s>\delta$,
$$
|u(t,x)-u(s,y)| \le c(\gamma,\beta,\delta,R, T,\varphi)\left(|(-A)^{\gamma}\left(x-y)\right)| + |t-s|^{\beta}\right)
$$
\item[(iv)] For any $t\in [0,T]$, $u(t,\cdot)\in{\mathcal E}$.
\item[(v)] $u(0,\cdot)=\varphi$.
\end{enumerate}
\end{Lemma}
{\bf Proof}. Let $R>0$, $\delta >0$ and $t,s\in [\delta,T]$, $x,y\in K_R.$ Then by Proposition \ref{p4.1}
and \ref{pc}  it follows that, for $\beta < \min\{g/2,1/2\}$,
\begin{equation}\label{3.4}
\begin{array}{l}
\ds{|u_m(t,x)-u_m(s,y)| \le c(\delta, T) \|\varphi\|_{\mathcal E}
\left(|(-A)^{\gamma}\left(x-y)\right)| + |t-s|^{\beta}\right)
}
\end{array}
\end{equation}
From the  the Ascoli--Arzel\`a theorem and
a diagonal extraction
argument, we deduce that
we can construct a subsequence such that
$$
u_{m_k}(t,x)\to u(t,x),
$$
uniformly in $[\delta,T]\times K_R$ for any $\delta>0,R>0.$ So that (i) follows.
Moreover, taking the limit in \eqref{3.4}, we deduce (iii).
We define $u(0,\cdot)=\varphi$. We can take the limit  $m_k\to \infty$ in Proposition \ref{pc}
with $t_1=t>0$ and $t_2=0$
and obtain for $x\in K_R$
$$
|u(t,x)-\varphi(x)|\le c(\varphi,R)\left( t^\beta + \left|A(e^{tA}-I)x\right|\right)
$$
This proves (ii) thanks to the strong continuity of
$(e^{At})_{t\ge 0}$. Finally, (iv) is an obvious consequence of
Proposition \ref{p4.1}.
\hfill $\Box$

\begin{Remark}
We could prove that $u$ is differentiable. This follows from a priori estimate on the
modulus of continuity of $Du_m$ and Ascoli--Arzel\`a
theorem. This further a priori estimate is rather technical but its proof does
not require new ideas.
\end{Remark}
Assume that $\varphi\in C^1_b(D(A);\R)$, the limit $u$ of the
subsequence $(u_{m_k})$ constructed above   may depend on the
choice of $(m_k)$. Therefore, at this point it is not clear that
this is any helpful to  construct a transition semigroup
$P_t\varphi$ for $\varphi\in B_b(D(A);\R)$. To avoid this problem, we shall
use the fact that we know the existence of a martingale solution
and of a stationary solution of \eqref{e1.2}. The stationary
solution will provide a candidate for the invariant measure $P_t$.

In order to emphasize the dependence on the initial datum, we shall denote by $u_m^\varphi$ the solution of
\eqref{e1.5}.

To prove Theorem \ref{t7.1} we need  further a priori estimates on
$X_m(t,x)$. The following Lemma is proved exactly as Lemma 7.4 in \cite{DPD-NS3D}.
 \begin{Lemma}
\label{l7.0}
For any $\delta\in (\frac 12\;,1+g]$, there exists a constant $c(\delta)>0$ such that for any $x\in H$, $m\in \N,$
and $t\in [0,T]$:
\begin{enumerate}
\item[(i)] $\ds{\E(|X_m(t,x)|^2)+\E\int_0^t |(-A)^{1/2}X_m(s,x)|^2ds\le |x|^2+t\;\mbox{\rm Tr}Q.
   }$

\item[(ii)] $\ds{\E\int_0^T\frac{|(-A)^{\frac{\delta+1}{2}}X_m(s,x)|^2}{(1+|(-A)^{\frac{\delta}{2}}X_m(s,x)|^2)
^{\gamma_\delta}}\;ds\le c(\delta),
} $ with $ \gamma_\delta=\frac{2}{2\delta-1}$ if $\delta\le 1$ and
$ \gamma_\delta=\frac{2\delta+1}{2\delta-1}$ if $\delta>1.$
\end{enumerate}
\end{Lemma}
\noindent It is well known that Lemma \ref{l7.0} can be used to
prove that the family of laws $(\mathcal L(X_m(\cdot,x))_{m\in
\N}$ is tight in $L^2(0,T;D((-A)^{s/2}))$ for $s<1$ and in
$C([0,T];D((-A)^{-\alpha}))$ for $\alpha>0.$ Thus, by the
Prokhorov theorem, it has a weakly convergent subsequence
$(\mathcal L(X_{m_k}(\cdot,x))_{k\in \N}$. We denote by $\nu_x$
its limit. By the Skohorod theorem there exists a stochastic
process $X(\cdot, x)$ on a probability space $(\Omega_x,\mathcal
F_x,\P_x)$  which belongs to $L^2(0,T;D((-A)^{s/2}))$ for $s<1$
and in $C([0,T];D((-A)^{-\alpha}))$ for $\alpha>0,$ satisfying
\eqref{e1.2} and such that for any $x\in D(A)$
\begin{equation}
\label{estar} X_{m_k}(\cdot,x)\to X(\cdot,x),\quad \P_x\mbox{  a.s.},
\end{equation}
in $L^2(0,T;D((-A)^{s/2}))$ and in
 $C([0,T];D((-A)^{-\alpha})).$

Note that it is not straightforward to build a family of
solution $(X(t,x),(\Omega_x,\mathcal F_x,\P_x))_x$ which is Markov.
Indeed, the sequence $(m_k)_{k\in \N}$ in \eqref{estar}
 may depend on $x$.
The key idea is to use the following Lemma which states that  the
sequence $(m_k)_k$ in Lemma \ref{l5.1} can be chosen independently
of $\varphi$.
\begin{Lemma}
\label{l7.1} There exists a sequence $(m_k)_{k\in \N}$ such that
for any $\varphi$ in ${\mathcal E}$ we have
$$
u^\varphi_{m_k}(t,x)\to u^\varphi(t,x),\quad \mbox{\it uniformly in}\;[\delta,T]\times K_R\;  \mbox{\it for
any}\;\delta>0,R>0.
$$
\end{Lemma}
{\bf Proof.} The proof is the same as for Lemma 7.5 in \cite{DPD-NS3D}, we
indicate the main ideas. Let $D$ be a dense countable set of $D(A)$.  It follows
from a diagonal extraction argument that there exists a sequence
$(m_k)_k$ such that \eqref{estar} holds for any $x\in D$. Then, it
follows from Lemma \ref{l5.1}-i) that, for any $\varphi$ in
${\mathcal E}$ and any subsequence of $(u_{m_k}^\varphi)$, we can
extract a subsequence which converges to a continuous map
$u^\varphi$. Moreover, it follows from Lemma \ref{l7.0}-ii)  that
$$
 X_{m_k}(\cdot,x)\to X(\cdot,x) \quad \mbox{in } D(A),\; d\P\times dt\textrm{ a.s.}, \textrm{ for any } x\in D.
 $$
It follows that
$\E_x[\varphi(X(t,x))]$ is defined $dt$ a.s. for $x\in D$
and taking the limit in \eqref{e1.6}
\begin{equation}
\label{es11}
 u^\varphi(t,x)=\E_x[\varphi(X(t,x))] \quad dt\textrm{ a.s.}, \textrm{ for any } x\in D.
\end{equation}
Therefore, any two accumulation points of $(u_{m_k}^\varphi)$ coincide on $D$.
By continuity, there is only one and the whole sequence $(u_{m_k}^\varphi)$ converges
to
 $u^\varphi$. This ends the
proof.\hfill$\Box$

We   now fix the sequence  $(m_k)_{k\in \N}$ and define for $\varphi\in {\mathcal E}$:
$$
P_t\varphi(x)=u^\varphi(t,x),\quad t\in [0,T],\;x\in D(A).
$$
As in \cite{DPD-NS3D}, it is easily deduced that $ P_t^*\delta_x$
defines a unique probability measure on $D(A)$. Therefore we can define
$P_t\varphi(x)$ for $\varphi \in B_b(D(A);\R)$.

 Then, for any  $x\in D(A)$, we build a martingale solution $X(\cdot,x)$ by
extracting a subsequence $(m_k^x)_k$ of $(m_k)_k$ such that
\eqref{estar} holds. It follows that
\begin{equation}
\label{es2}  P_t\varphi(x)=\E_x[\varphi(X(t,x))],\quad x\in
D(A),\;t\in [0,T],
\end{equation}
provided $\varphi\in C_b(D((-A)^{-\alpha});\R)\cap \mathcal E$.
 It  easily checked that \eqref{es2} remains true for
$\varphi$ uniformly continuous in $D((-A)^{-\alpha})$. Thus $
P_t^*\delta_x$ - seen as a probability measure on
$D((-A)^{-\alpha})$ - is the law of $X(t,x)$. Since $
P_t^*\delta_x$ is a probability measure on $D(A)$,  i) of
Definition \ref{Def} follows. Moreover \eqref{es2} remains true
for $\varphi\in B_b(D(A);\R)$.

Recall that for any subsequence of $(m_k)_k$ and any $x\in D(A)$,
we have a subsequence $(m_k^x)_k$ and a martingale solution
$X(\cdot,x)$ such that $X_{m_k^x}(\cdot,x)\to X(\cdot,x)$ in law
in $C(0,T;D((-A)^{-\alpha}))$.

 To end the proof of Theorem \ref{t7.1}, we prove below the
 following result.
 \begin{Lemma}\label{lCyr}
Let $X(\cdot,x)$ be a limit process of a subsequence of
$(X_{m_k}(\cdot,x))_{m_k}$.  Then, for any  $(n,N)\in\N^2$,
$t_1,\dots,t_n \ge 0$ and $(f_k)_{k=0}^n \in
C_c^\infty(P_N H)$ (i.e. $f_k(x)=f_k(P_Nx)$), we have
\begin{equation}
\label{Cyr1}
\begin{array}{r}
\E_x\left(f_0(X(0,x))f_1(X(t_1,x)),\dots, f_n(X(t_1+\dots+t_n,x))\right)\quad\quad\quad\quad\quad\quad\quad\quad\\
=f_0(x)P_{t_1}\bigg[f_1
P_{t_2}\bigg(f_2P_{t_3}(f_3\dots)\bigg)\bigg](x).
\end{array}
\end{equation}
\end{Lemma}
By classical arguments, uniqueness in law of the limit
process follows from \eqref{Cyr1}. Combining uniqueness of the
limit and compactness of the sequence $(X_{m_k}(\cdot,x))_k$, we
obtain that $X_{m_k}(\cdot,x)\to X(\cdot,x)$ in law in
$C(0,T;D((-A)^{-\alpha}))$. It follows that the map $x\to\mathcal
P_x$ defined in Definition \ref{Def} depends measurably on $x$,
and that $\mathcal P_{X(t,x)}(\mathcal A)$ is a random
variable for any $\mathcal A$ as in Definition \ref{Def} (ii).

 We set, for $t_1,\dots,t_k,s_1,\dots,s_n\ge 0$ and $A_0,\dots,A_k,A'_0,\dots,A'_n\in {\mathcal B}(D(A))$,
$$
\Espace \left\{
\begin{array}{rcl}
\mathcal A&=&\{X(0)\in A_0,\dots,X(t_1+\dots+t_k)\in A_k\},\\
 \mathcal A'&=&\{X(0)\in A_0',\dots,X(s_1+\dots+s_n)\in A_n'\}.
\end{array}
\right.
$$
Since \eqref{Cyr1} is easily extended to any $(f_0,\dots,f_n)\in
(B_b(D(A);\R))^{n+1}$, we deduce
$$
\mathcal P_x\left(\mathcal A
\right)=1_{A_0}(x)P_{t_1}\bigg[1_{A_1}
P_{t_2}\bigg(1_{A_2}P_{t_3}(1_{A_3}\dots)\bigg)\bigg](x).
$$
Applying successively \eqref{Cyr1} to
$(1_{A_0'},\dots,1_{A_n'},x\mapsto\mathcal
P_x\left(\mathcal A' \right))$, $s_1,\dots,s_n, t-(s_1+\dots s_n)$ and to
$(1_{A_0'},\dots,1_{A_n'},1_{A_0},\dots,1_{A_k})$, $s_1,\dots,s_n,t-(s_1+\dots s_n), t_1,\dots,t_n$,
we obtain
\begin{equation}\label{Cyr4}
\P_x\left(X(\cdot,x)\in \mathcal A',X(t+\cdot,x)\in \mathcal
A\right) =\E_x\left(1_{\mathcal A'}(X(\cdot,x))\mathcal
P_{X(t,x)}(\mathcal A)\right).
\end{equation}
provided $t\geq s_1+\dots+s_n$.
This yields
(ii) of Definition \ref{Def}.

{\bf Proof of Lemma \ref{lCyr}.} We set $f(x_0,\dots,x_n)=f_0(x_0)\dots f_n(x_n)$. It follows from the weak Markov
property that
\begin{equation}\label{Cyr7}
\E_x^{m_k^x}\left(f(X_{m_k^x}(0,x),\dots,X_{m_k^x}(t_1+\dots+t_n,x))\right)=
f_0(x)P_{t_1}^{m_k^x}\bigg[f_1
P_{t_2}^{m_k^x}\bigg(f_2P_{t_3}^{m_k^x}(f_3\dots)\bigg)\bigg](x).
\end{equation}
Moreover, the convergence in $C(0,T;D((-A)^{-\alpha}))$ gives
\begin{equation}\label{Cyr6}
\E_x\left(f(X(0,x),\dots, X(t_1+\dots+t_n,x))\right)=\lim
\E_x^{m_k^x}\left(f(X_{m_k^x}(0,x),\dots)\right).
\end{equation}
Therefore it remains to prove that
\begin{equation}\label{Cyr5}
f_0(x)P_{t_1}^{m_k}\bigg[f_1
P_{t_2}^{m_k}\bigg(f_2P_{t_3}^{m_k}(f_3\dots)\bigg)\bigg](x) \to
f_0(x)P_{t_1}\bigg[f_1
P_{t_2}\bigg(f_2P_{t_3}\dots\bigg)\bigg](x),
\end{equation}
uniformly on $D(A)$.

We prove \eqref{Cyr5} by induction on $n\in\N$. For $n=0$, this is
 trivial. For $n=1$, this is Lemma \ref{l7.1}. Assume that \eqref{Cyr5}
 is true for $n\in\N$. We set
 $$
\Espace \left\{
\begin{array}{rcl}
I_{m}(x)&=&f_0(x)P_{t_1}^{m}\bigg[f_1
P_{t_2}^{m}\bigg(f_2P_{t_3}^{m}(f_3\dots)\bigg)\bigg](x) -
f_0(x)P_{t_1}\bigg[f_1
P_{t_2}\bigg(f_2\dots\bigg)\bigg](x),\\
J_{m}(x)&=&f_1(x)P_{t_2}^{m}\bigg[f_2
P_{t_3}^{m}\bigg(f_3P_{t_4}^{m}(f_4\dots)\bigg)\bigg](x) -
f_1(x)P_{t_2}\bigg[f_2
P_{t_3}\bigg(f_3\dots\bigg)\bigg](x).\end{array} \right.
 $$
Remark that
\begin{equation}\label{Cyr8}
I_{m_k}(x)=f_0(x)\left(
P_{t_1}^{m_k}J_{m_k}(x)+(P_{t_1}^{m_k}-P_{t_1})\bigg[f_1
P_{t_2}\bigg(f_2P_{t_3}\dots\bigg)\bigg](x))\right).
\end{equation}
By the induction argument, $J_{m_k}(x)\to 0$ uniformly. Hence
$$
||f_0 P_{t_1}^{m_k}J_{m_k}||_0\leq ||f_0||_0 ||J_{m_k}||_0\to 0.
$$
 Moreover, since $f_0=0$ out of a bounded set,
then it follows from Lemma \ref{l7.1}
$$
f_0(x)(P_{t_1}^{m_k}-P_{t_1})\bigg[f_1
P_{t_2}\bigg(f_2P_{t_3}\dots\bigg)\bigg](x)\to 0,
$$
uniformly on $D(A)$, which yields \eqref{Cyr8}. \hfill $\Box$

\section{Proof of Theorem \ref{t7.2}}

Now we observe that, since the noise is nondegenerate, then
$P^{m_k}_t$ has a unique invariant measure $\nu_{m_k}.$ Moreover,
we have the following result proved as in \cite{DPD-NS3D} (see
Lemma 7.6).
\begin{Lemma}
\label{l6.2} There exists a constant $C_1$ such that for any $k\in
\N$
$$
\int_H\left[|(-A)^{1/2}x|^2+|Ax|^{2/3}+|(-A)^{1+\frac{g}{2}}x|^{\frac{1+2g}{10+8g}}\right]\nu_{m_k}(dx)<
C_1.
$$
\end{Lemma}
It follows that the sequence $(\nu_{m_k})_{k\in \N}$ is tight on
$D(A)$ and there exists a subsequence, which we still denote by
$(\nu_{m_k})_{k\in \N}$, and a measure $\nu$ on $D(A)$ such that
$\nu_{m_k}$ converges weakly to $\nu.$ Moreover
$\nu(D((-A)^{1+\frac{g}{2}})=1$.

Let us take $\varphi\in {\mathcal E}$.  It follows from the
invariance of $\nu_m$ and the convergence properties of the
approximations of $P_t\varphi$ that for any $t\ge 0$
\begin{equation}
\label{invariance} \int_H P_t\varphi(x)\nu(dx)=\int_H
\varphi(x)\nu(dx).
\end{equation}
Therefore $\nu$ is an invariant measure.
The strong Feller property is a consequence of Proposition
\ref{p4.1}.  Hence, by Doob Theorem, the strong mixing property is
a consequence of the irreducibility.
 This latter property is implied by the following Lemma. Its
proof significantly differs from the additive case treated in
\cite{DPD-NS3D}.
\begin{Lemma}
\label{irreducibility}
Let $x_0\in D(A)$, $\varepsilon>0$
and $\varphi\in {\mathcal E}$ be such that
$\varphi(x) =1$ for $x$ in $B_{D(A)}(x_0,\varepsilon)$, the
ball in $D(A)$ of center $x_0$ and radius $\varepsilon$. Then
for any $t>0$ and $x\in D(A)$ we have $P_t\varphi(x)>0$.
\end{Lemma}
{\bf Proof.}  It is classical that, since $x\in D(A),$ there
exists $T^*>0$ and $\overline{x}\in C([0,T^*];D(A))$ such that
$$
\overline{x}(t)=e^{tA}x+\int_0^te^{(t-s)A}b(\overline{x}(s))ds, \quad t\in [0,T^*].
$$
This follows from a fixed point argument. Moreover, it is not difficult to see that
$\overline{x}\in  L^2(0,T^*;D((-A)^{3/2}))$, so that $\overline{x}(t)\in D((-A)^{3/2})$ a.e. and we may change $T^*$
so that $\overline{x}(T^*)\in D((-A)^{3/2})$. Morever, we can assume that $x_0\in D((-A)^{3/2})$.
Then we set $\overline{g}=0$ on $[0,T^*]$ and define $\overline{x}$, $\overline{g}$
 on $[T^*,T]$ by
$$
\overline{x}(t)=\frac{T-t}{T-T^*}\;\overline{x}(T^*)+\frac{t-T^*}{T-T^*}\;x_0,\quad \overline{g}(t)=\frac{d\overline{x}}{dt}\;-A\overline{x}-b(\overline{x}).
$$
We also define
$$
R=2\sup_{t\in [0,T]} |A\overline{x}(t)|,\quad
b_R(x)=\vartheta\left(\frac{|Ax|}R\right)b(x),
$$
where $\vartheta \in C^\infty_0(\R)$ takes its values in $[0,1]$,
is equal to $1$ on $[0,1]$ and vanishes on $[2,\infty)$.
 Then, for
$n\in \N$, we write
$$
\Delta t = T/n,\; t_k=k\Delta t, \mbox{ for } k\in \N, \mbox{ and } \dot {W}_{n}(t) =\frac{ P_n  W(t_k)
- P_n W(t_{k-1})}{\Delta t},\quad
t\in [t_k,t_{k+1}).
$$
It is not difficult to see that the equation
\begin{equation}
\label{3.6}
\left\{
\begin{array}{lr}
\ds{dX^{n,R}}\ds{= (AX^{n,R} +b_R(X^{n,R})+P_n \overline{g} }&\ds{-\Phi_n(e^{A(t-t_{k-1})}X^{n,R}(t_{k-1}))\dot{W}_n)dt } \\
&\ds{ +\Phi(X^{n,R})dW,}\; t\in [t_k,t_{k+1}),\\
X^{n,R}(0)=x&,
\end{array}\right.
\end{equation}
has a unique solution $X^{n,R}$  in $C([0,T];D(A))$. We set
$$
G^n(t)=\Phi(X^{n,R}(t))^{-1}\left( P_n \overline{g}(t)-\Phi_n(e^{A(t-t_{k-1})}X^{n,R}(t_{k-1}))\dot{W}_n(t)\right), \; t\in [t_k,t_{k+1}).
$$
By Girsanov Theorem
$\widetilde{W}_n(t) = W(t) -\int_0^t G^n(s)(s) ds$
defines a cylindrical Wiener process for the probability measure
$$
d\P_n=\exp \left[\int_0^T\left( G^n(s),dW(s)\right)-\frac{1}{2}\int_0^T\left| G^n(s)\right|^2ds\right]d\P.
$$
It follows that the laws of $X^{n,R}$ and $X^R$ are equivalent, where $X^R$ is the solution
of
\begin{equation}
\label{3.7}
\left\{
\begin{array}{l}
\ds{dX^{R}= (AX^{R} +b_R(X^{R}))dt +\Phi(X^{R})dW, \;
t\in [t_k,t_{k+1}),}\\ 
X^{R}(0)=x.
\end{array}\right.
\end{equation}
Moreover, we have, for a subsequence $m_k$,
$$
\begin{array}{ll}
P_t\varphi(x)&\ds{=\lim_{m_k\to \infty} \E(\varphi(X_{m_k}(t,x)))}\\
&\ds{\ge \lim_{m_k\to \infty} \E(\varphi(X_{m_k}(t,x))1\hspace{-0.1cm}\mbox{I}_{\tau_{m_k}^R\ge T})}\\
&\ds{\ge \lim_{m_k\to \infty} \E(\varphi(X_{m_k}^R(t,x))1\hspace{-0.1cm}\mbox{I}_{\tilde\tau_{m_k}^R\ge T})}
\end{array}
$$
where $\tau_{m_k}^R=\inf\{t\in [0,T], \; |AX_{m_k}(t,x)|\ge R\}$,
$\tilde\tau_{m_k}^R=\inf\{t\in [0,T], \; |AX_{m_k}^R(t,x)|\ge R\}$, and $X_{m_k}^R$ is the solution of
the Galerkin approximation where $b$ has been replaced by $b_R$.
Since, it easy to check that
$$
\lim_{m_k\to \infty} \E(\varphi(X_{m_k}^R(t,x))1\hspace{-0.1cm}\mbox{I}_{\tilde\tau_{m_k}^R\ge T})
= \E(\varphi(X^R(t,x))1\hspace{-0.1cm}\mbox{I}_{\tau^R\ge T})
$$
where $\tau^R=\inf\{t\in [0,T], \; |AX^R(t,x)|\ge R\}$, we deduce that it is sufficient to prove that
$
\E(\varphi(X^R(t,x))1\hspace{-0.1cm}\mbox{I}_{\tau^R\ge T})>0
$
and, since the laws of $X^{n,R}$ and $X^R$ are equivalent, that
$
\E(\varphi(X^{n,R}(t,x))1\hspace{-0.1cm}\mbox{I}_{\tau^{n,R}\ge T})>0
$
with $\tau^{n,R}=\inf\{t\in [0,T], \; |AX^{n,R}(t,x)|\ge R\}$.
We prove below that $X^{n,R}$ converges to $\overline{x}$ in $L^2(\Omega;C([0,T];D(A)))$
as $n\to \infty$. Since
$\varphi(\overline{x}(T))=1$, the claim follows.

To prove that $X^{n,R}$ converges to $\overline{x}$, we first observe that
 $X^{n,R}$ is uniformly bounded with respect to $n$ in $L^2(\Omega;C([0,T];D((-A)^\gamma)))$
for $\gamma <1+g/2$. This is proved thanks to the integral form of \eqref{3.6} and similar argument
as in Proposition \ref{P2.2}. This enables to control the difference between
$ \Phi_n(e^{A(t-t_{k-1})}X^{n,R}(t_{k-1}))$ and $\Phi(X^{n,R}(t))$, for $t\in [t_k,t_{k+1}]$.

Then, we write the difference between the integral equation satisfied by
$X^{n,R}$ and $\overline{x}$ and see that it is sufficient to prove that $\int_0^t e^{A(t-s)} \Phi_n(e^{A(s-s_{k-1})}X^{n,R}(s_{k-1}))\dot{ W}_n ds
-\int_0^t e^{A(t-s)} \Phi(X^{n,R}(s)dW(s)$ goes to  zero in $L^2(\Omega;C([0,T];D(A)))$, where
$s_{k_1}=t_{k-1}$ for $s\in [t_k,t_{k+1}]$. This latter
point is not difficult to prove.
\hfill $\Box$

Let $(X_{m_k}(\cdot,\nu_{m_k}))_k$ be the  sequence of  stationary
solutions. Proceeding as in the end of section $4$, we remark that
to prove the convergence, it is sufficient to establish uniqueness
in law of the limit of subsequence of
$(X_{m_k}(\cdot,\nu_{m_k}))_k$. Remark that Theorem
\ref{t7.2}-(ii) implies such uniqueness.

So, to conclude Theorem \ref{t7.2}, it remains to establish (ii).
By classical arguments, it follows from Lemma \ref{lCyr} that it
is sufficient to establish that
\begin{equation}\label{Cyr9}
\begin{array}{r}
\E_{\nu}\left(f_0(X(0,\nu))\dots
f_n(X(t_1+\dots+t_n,\nu))\right)\quad\quad\quad\quad\quad\quad\quad\quad\quad\quad\\
=\int_{D(A)}f_0(x)P_{t_1}\bigg[f_1
P_{t_2}\bigg(f_2P_{t_3}(f_3\dots)\bigg)\bigg](x)\,\nu(dx).
\end{array}
\end{equation}
 for any $(t_1,\dots,t_n)$ and $(f_0,\dots,f_n)$  as in Lemma
\ref{lCyr}.

Remark that
$$
\begin{array}{r}
\E_{\nu_{m_k}}\left(f_0(X_{m_k}(0,\nu_{m_k}))\dots
f_n(X_{m_k}(t_1+\dots+t_n,\nu_{m_k}))\right)\quad\quad\quad\quad\quad\quad\quad\quad\quad\quad\\
=\int_{D(A)}f_0(x)P_{t_1}^{m_k}\bigg[f_1
P_{t_2}^{m_k}\bigg(f_2P_{t_3}^{m_k}(f_3\dots)\bigg)\bigg](x)\,\nu_{m_k}(dx).
\end{array}
$$
The convergence of the right-hand-side comes from the convergence
in law in \linebreak$C(0,T;D((-A)^{-\alpha}))$. Applying
\eqref{Cyr5}, we can take the limit in the left-hand-side and then
conclude.

\end{document}